\documentclass[a4paper,12pt]{article} 
\usepackage{amsmath,graphicx,amssymb,fancyhdr,amsthm,verbatim}

\theoremstyle{definition}  
  
\theoremstyle{remark}  
  
\def\beq{\begin{eqnarray}}  
\def\eeq{\end{eqnarray}}  
\def\bsp{\begin{split}}  
\def\esp{\end{split}}

\begin{document}   
   
\title{\Large\textbf{Basis for scalar curvature invariants in three dimensions}}   
\author{{\large\textbf{A. A. Coley$^{1}$ A. MacDougall$^{1}$  and D. D. McNutt$^{1}$ } }  
 \vspace{0.3cm} \\   
$^{1}$Department of Mathematics and Statistics,\\   
Dalhousie University,    
Halifax, Nova Scotia,\\    
Canada B3H 3J5    
\vspace{0.2cm}\\ 
\texttt{aac@mathstat.dal.ca, amacdougall,ddmcnutt@dal.ca} }   
\date{\today}   
\maketitle   
\pagestyle{fancy}   
\fancyhead{} 
\fancyhead[EC]{A. Coley and A. MacDougall}   
\fancyhead[EL,OR]{\thepage}   
\fancyhead[OC]{Scalar Curvature Invariants}   
\fancyfoot{} 
   
\begin{abstract}

$\mathcal{I}$-non-degenerate spaces 
are spacetimes that can be characterized uniquely by their scalar curvature invariants.
The ultimate goal of the current work is to construct a basis for the scalar polynomial curvature invariants in 
three dimensional Lorentzian spacetimes. In particular, we seek a minimal
set of algebraically independent scalar curvature invariants formed by the contraction of
the Riemann tensor and its covariant derivatives up to fifth order of differentiation. 
We use
the computer software \emph{Invar}  to calculate an overdetermined basis
of scalar curvature invariants in three dimensions.
We also discuss the
equivalence method and the Karlhede algorithm for computing Cartan invariants
in three dimensions.

\end{abstract}   

\section{Introduction}

Scalar curvature invariants are scalars constructed from the Riemann tensor and its covariant
derivatives.  Scalar curvature invariants can be used to study the inequivalence of metrics and
curvature singularities \cite{exactSol}.  In the case of the Lorentzian $\mathcal{I}$-non-degenerate
spaces, they can be used to determine the equivalence of metrics \cite{inv}.
Scalar curvature invariants have been utilised in the study of universal spacetimes \cite{CGHP}, and in VSI and CSI spacetimes
\cite{coley,CSI,applicat}.  
In particular, scalar curvature invariants have been studied due to their potential use in general
relativity \cite{exactSol}.  Primarily the scalar
curvature invariants formed from the Riemann tensor $R_{abcd}$ only (contractions involving products of
the undifferentiated Riemann tensor only, the so called \emph{algebraic invariants}) has been investigated in
four-dimensional (4D) Lorentzian spacetimes.

Scalar curvature invariants are of primary importance in $\mathcal{I}$-non-degenerate spacetimes
\cite{inv,coley,highSI}. Indeed, these
spacetimes can be completely characterized by their scalar curvature invariants.  This leads to the
natural problem of attempting to find a basis for the scalar curvature invariants formed from the Riemann tensor
(up to some order of covariant differentiation).  Much work has gone into the problem of constructing a
basis in the 4D Lorentzian algebraic case.  In this case one can form 14 functionally independent
scalar curvature invariants.  The smallest set that contains a maximal set of algebraically independent
scalars consists of $17$ polynomials \cite{zak}.  As reviewed in the introduction of \cite{roman}, a
number of proposed independent sets were given by Narlikar and Karmarkar \cite{narlikar}, Geheniau and
Debever \cite{geheniau}, Petrov \cite{petrov}, and others.  All of these sets were shown to be
deficient for various reasons. A set of algebraic invariants was presented by
Carminati and McLenaghan \cite{roman}, consisting of 16 curvature invariants, that contains invariants of lowest
possible degree and contains a minimal set for any Petrov type and for any specific choice of Ricci
tensor type in the perfect fluid and Einstein-Maxwell cases.  In general, the expressions relating
invariants to the basis members of an independent set can be very complicated, and can be singular in
certain algebraic cases.

Much less work has gone into studying curvature 
invariants formed using covariant derivatives of the Riemann tensor (differential invariants)
\cite{CSI,FKWC1992,DecaniniFolacci2007,invar,diffInvar}.  Our ultimate goal is to provide a basis of 
all such scalar curvature invariants.
This is primarily a mathematical question.  However, it also has applications in mathematical physics.

There are differing notions of what is meant by a basis.  For example, the number $N(n,p)$ of algebraically independent
quantities formed from the first $p$ derivatives of the Riemann tensor in dimension $n$ can be determined (see equation
(\ref{thomas}) below).  These $N(n,p)$ components correspond to the independent components of the Riemann tensor and
its covariant derivatives.  Hence, in principal, a basis could consist of these $N(n,p)$ scalars because all other
curvature invariants are functions of these scalars.  However, the expressions relating these invariants can be very
complicated, involving roots of high order, and they may therefore be singular in certain algebraic cases.  For actual
classification using invariants, a different type of basis is needed.  The common solution to this problem is to seek a
basis of scalars such that all other scalars are polynomials in this basis.  This is the approach that has been widely
used in the study of scalar curvature invariants in 4D Lorentzian signature.

Rather than searching for an independent basis, we may search for a \emph{complete} basis
$\{I_1,I_2,\ldots,I_n\}$, meaning that any scalar curvature invariant can be expressed as a polynomial
of the elements of $\{I_1,I_2,\ldots,I_n\}$ but no member of $\{I_1,I_2,\ldots,I_n\}$ can be expressed
in terms of the others.  For example, the smallest known complete set in the 4D
algebraic Lorentzian case consists of $38$ scalars \cite{sneddon}.

The main motivation for studying this problem in three dimensions (3D)
is mathematical, although there are some possible applications to
mathematical physics \cite{AlievNutku,Hall}.  In particular, it may be applied to the initial
value problem in general relativity; in this case a unique solution is
given by a 3-metric and 3D extrinsic curvature specified on a Cauchy
hypersurface.  Hence, an important step is to characterize the
hypersurface using scalar curvature invariants.

The primary motivation for studying scalar curvature
invariants comes from the $\mathcal{I}$-non-degenerate theorem.

\subsection{$\mathcal{I}$-non-degenerate spacetimes}  
In \cite{inv}  the class of 4D Lorentzian
manifolds that can be completely characterized by the scalar polynomial
curvature invariants constructed from the Riemann tensor and its covariant derivatives was determined.
The notion of an $\mathcal{I}$-non-degenerate spacetime metric was introduced, which implies that the spacetime metric is locally determined
by its curvature invariants. It was proven  that a spacetime
metric is either $\mathcal{I}$-non-degenerate or a degenerate Kundt metric. Therefore, a metric
that is not characterized by its curvature invariants must be of Kundt
form. This theorem is also true in 3D, and it is also likely true in arbitrary dimenions;  
in \cite{highSI}, a number of results that generalize these results to higher dimensions were introduced, and their
consequences and potential physical applications were discussed.

Let us descibe this theorem in a little more detail.
A (one-parameter) \emph{metric deformation} $\hat{g}_\tau, \tau\in[0,\epsilon)$ of a 
spacetime $(\mathcal{M},g)$ is a family of smooth metrics on $\mathcal{M}$ such that
(i) $\hat{g}_\tau$ is continuous in $\tau$,
(ii) $\hat{g}_0 = g$,
(iii) $\hat{g}_\tau$ for $\tau>0$ is not diffeomorphic to $g$.
Given a spacetime $(\mathcal{M},g)$ with a set of invariants $\mathcal{I}$, then if there does not exist a metric deformation of $g$ with the same set of invariants as $g$, then we will call the set of invariants \emph{non-degenerate}. The spacetime metric $g$ will be called \emph{$\mathcal{I}$-non-degenerate.}
In \cite{inv}, it was proven that a spacetime metric 
{\em{is either $\mathcal{I}$-non-degenerate or the metric is a Kundt metric}}; that is,
a Lorentzian metric that is $\mathcal{I}$-non-degenerate is locally characterized uniquely by its 
scalar curvature invariants.

\newpage

\section{Bounds on number of covariant derivatives}

We are interested in constructing a basis for the scalar curvature invariants up to some order of covariant differentiation in general
Lorentzian spacetimes (although in this paper we are specifically interested in 3D).  In other words, we seek a minimal (or
overdetermined) set of algebraically independent scalar curvature invariants formed by the contraction of the Riemann tensor and its
covariant derivatives up to some order of differentiation $p$.  Hence, we must determine the bound on $p$.  Obviously, we seek the
minimum value of $p$ such that all geometrical information of the spacetime can be obtained.


The bound  $q \equiv p+1$ is used in the Cartan-Karlhede algorithm for determining the equivalence of spacetimes \cite{karlhede} . The bound on the algorithm (by considering the `worst-case scenario' for the number of steps required for the algorithm to finish) is 
\begin{equation}q  = N_0 + n + 1  \label{bound}\end{equation}
 where $N_0$ is the dimension of the isotropy (automorphism)
 group of the Riemann tensor and $n$ is the dimension of the spacetime.
At the final step $q$ of the Cartan-Karlhede algorithm, new classification functions 
can be introduced.
Therefore, in the classification
we need the Riemann tensor and its first $q$
covariant derivatives.

\subsection{Bounds in 4D}  
We first recall some results in 4D.  Cartan's
argument gives a bound of $n(n+1)/2 = 10$ \cite{exactSol}.  It is, however, impossible in the 4D
Lorentzian case to have $p>7$ (except for the constant curvature case).  To see this suppose that $p>7$.  In this case, in
the Karlhede algorithm we must find at most ten functions on the ten dimensional frame bundle (6 Lorentz group parameters and
four independent functions of spacetime).  Now, since we need $8$ or $9$ steps to terminate the algorithm, then at most two
of the parameters were already fixed at the beginning.  Then the undetermined part of the Lorentz group would be of dimension
at least $4$, but there is no choice of Weyl and Ricci curvature allowing this (i.e.,  invariant under a 4D subgroup of the
Lorentz group).  Then the curvature would be invariant under the whole Lorentz group, which implies that the spacetime is of
constant curvature \cite{exactSol}.  Furthermore, in 4D we know that if $p\geq 6$ then the spacetime is degenerate Kundt
\cite {highSI}.  Therefore, for $\mathcal{I}$-non-degenerate spacetimes, we have $p\leq 5$
($q\leq 6)$.  We note that $N(4,6) = 2094$
(see equation (\ref{thomas}) below).

It is conjectured that for $\mathcal{I}$-non-degenerate spacetimes all
Cartan invariants are determined (up to possible discrete complex
transformations) by scalar polynomial curvature invariants (as in the
Riemannian case) \cite{inv}.

In 4D Lorentzian signature, a complete basis at zeroth order has 38 objects and
contains Riemann scalars up to degree eleven \cite{sneddon}.  In
\cite{invar}, \emph{Invar} was used to determine 25 of the 38 objects
in this basis.  In this basis there are 24 syzygies between the
invariants.  The syzygies are polynomials, but they are not yet known.
If \emph{Invar} were able to be run to an arbitrary degree they would
in principal be found, but this is not computationally feasible. Recently,
using techniques from graph theory, Carminati and Lim
claim to be able to reproduce Sneddon's basis \cite{zak}.

\subsection{Bounds in 3D}
The question of the bound for $p$ in 3D was addressed by Sousa, Fonseca, and Romero in \cite{Sousa}. In 3D spacetimes, the Weyl tensor vanishes, and the canonical frame of the Karlhede algorithm \cite{karlhede} is aligned with principal directions of the Ricci tensor rather than the Weyl tensor. 
All spacetimes have $N_0\leq 1$, hence we have $q=p+1\leq 5$ (the same as the 3D Riemannian case, which has no two 
dimensional isotropy group). See table \ref{3dTable} for the (dimension of the) isotropy groups of spacetimes of different 
Segre types of the Ricci tensor \cite{Sousa}, and we also show the related Ricci type \cite{CSI,hall}. 

It is common to discuss the Segre types of the traceless Ricci tensor \cite{exactSol} whence the
case $\{(1,11)\}$ with full 3D Lorentz subgroup corresponding to spacetimes of constant curvature is of
Ricci type O.  There is a 2D isotropy group of the Lorentz group spanned by a boost and a null rotation
(either of these together with a spatial rotation gives rise to the full 3D Lorentz group).  The 2D
subgroups arise as special cases of Segre types $\{(1,1),1\}$ and $\{(2,1)\}$ in table \ref{3dTable}
with 1D isotropy group (with an additional isotropy).
However, these resulting spacetimes must be degenerate Kundt (i.e., not $\mathcal{I}$-non-degenerate), 
and are not considered further here. Suppose that the Karlhede algorithm takes five steps to complete (i.e., $p=4$). 
Then at most two of the parameters on the fibre bundle were fixed at the beginning. Then the undetermined part of the Lorentz group has $N_0\geq 1$, and we see that this argument does not improve the bound in 3D as it does in the 4D case.

It was shown in \cite{MW2013} that
this bound is sharp (i.e., $q=5$ is attained); that is,  a 3D
Lorentzian manifold exists for which the fifth covariant derivative of
the Riemann tensor is required to classify the spacetime invariantly.
If the isotropy group is one dimensional, then $q$ can be 5. 
However, in  the (general) cases (see the table) that do not have a one dimensional
isotropy group (and are not degenerate Kundt), 
then  $q \leq 4$. Thus, in these cases the
classification is simpler. We note that $N(3,4) = 147$ and $N(3,5) = 228$.
We might be able to improve upon $q=4$ on a case by case basis. Indeed, there is a sense in which in general 
$q=1$; note that $N(3,1) = 18$.

\begin{table}
\caption{Ricci type to Segre type conversion scheme. }
\centering
\begin{tabular}{|c|c|c|c|c|}
\hline
Segre type &  & Ricci type & Isotropy group & Dimension of Isotropy Group \\
\hline\hline
$\{1,11\}$ & & $I$ & none & $0$\\
\hline
$\{1,(11)\}$ & & $I$ & $SO(2)$ & $1$\\
\hline
$\{(1,1)1\}$ & & $D$ & $SO(1,1)$ & $1$\\
\hline
$\{(1,11)\} $ & $\lambda_1\neq 0$ & $D$ & $SO(2,1)$&$3$ \\
              & $\lambda_1=0$ & $O$ & & \\
\hline
$\{1\bar{z}z\}$ && $I$ & none & $0$\\
 \hline
$\{21\}$ & & $II$ & none&$0$\\
\hline
$\{(21)\}$ & $\lambda_1\neq 0$& $II$ & null-rotations &$1$\\
           & $\lambda_1= 0$ & $N$ & & \\
\hline
$\{3\}$ & $\lambda_1\neq 0$& $II$ & none &$0$ \\
           & $\lambda_1= 0$ & $III$ & &\\
\hline
\end{tabular}\label{3dTable}
\end{table}

\subsection{Cartan invariants in 3D}

As noted earlier, it has been conjectured that for $\mathcal{I}$-non-degenerate spacetimes all
Cartan invariants are determined by scalar polynomial curvature invariants  \cite{inv}.
Hence the equivalence method suggests an approach to determine a basis of 
all scalar curvature invariants.

In particular, 
the question of the maximal order of covariant derivative required for the 
invariant classication of a pseudo-Riemannian manifold, as discussed above \cite{Sousa,MW2013}, is relevant
in determining the worst case scenarios for implementing the equivalence algorithm.
This then provides an upper bound on the number of covariant derivatives
needed for a basis of 
all scalar curvature invariants.

The equivalence problem in 3D was first considered in \cite{Sousa},  and
the bound $q=5$ was determined, and using methods of
two-component real spinors the isotropy groups and the canonical
forms of algebraic classifications of the Ricci and Cotton-York
spinors were given.  The results were applied to the equivalence of
Godel-type spacetimes in 3D \cite{Sousa} (also see \cite{McNutt}).
A further, more comprehensive, example will be presented
later.

\newpage

As an illustration, let us consider the 3D  (P-type I)
G\"{o}del-like metric as given in  \cite{Sousa}
	\begin{equation}
		ds^2 = - (dt + H(r) d\phi)^2 + D^2(r) d\phi^2 + dr^2.
	\end{equation}
The Cartan invariants consist of three zeroth order invariants, one first order invariant
(essentially the spin coefficient $\tau$), and two frame derivatives of the 
zeroth order invariants \cite{McNutt}: $\left\{I_0, I_1, I_2, {I_0}',  {I_1}', {I_2}'\right\}$, where
	\begin{align*}
		I_0 = \frac{H'}{D},~~ I_1 = \frac{D''}{D},~~I_2 = \frac{D'}{D}
	\end{align*}
and a prime denotes differentiation with repect to r.

The three zeroth order polynomial scalar curvature invariants 
expressed as polynomials 
in the Cartan invariants are \cite{McNutt}:
	\begin{align*}
		{R^a}_a
		\;		
		&=
		-2\,{I_1} +  \frac{1}{2} \, I_0 ^{2}
	\end{align*}
	\begin{align*}
		R^{ab} R_{ab}
		\;
		&=
		2\,  I_1  ^{2}-2\,I_1  I_0^{2}- \frac{1}{2}\, \left(I_0'  \right) ^{2}+ \frac{3}{4}\, I_0^{4}
	\end{align*}
	\begin{align*}
		R^{ab} {R_a^{~c}} R_{bc}
		\;
		&=
		-2\,  I_1  ^{3}+ \frac{3}{4}\,I_1  \left(I_0'  \right) ^{2}+3\, I_1  ^{2} I_0 ^{2}- \frac{3}{2}\,I_1  I_0^{4}+ \frac{1}{8}\,  I_0   ^{6}
	\end{align*}
The first order scalar invariants, expressed as polynomials in the Cartan invariants, 
can be written down explicitly:
	\begin{align*}
		{R^{;a}}_{;a}
		\;
		&=
		\left( -I_0 I_0' +2\,{I_1}'  \right) ^{2},
	\end{align*}
and where	
$R^{bc;d} R_{bc;d}$ and
$R^{bc;d} R_{bd;c}$ are more complicated explicit polynomial functions of 
$\left\{I_0, I_1, I_2, {I_0}',  {I_1}', {I_2}'\right\}$  \cite{McNutt}. Thus the 3D 
G\"{o}del spacetimes are characterized by these six zeroth and first order scalar curvature
invariants.

This example is intended to illustrate the applicability of the algebraic               
relations of the Cartan invariants at each order to the simplification and study of                  
polynomial scalar curvature invariants. In 3D these curvature invariants are           
expressed in terms of the Ricci components; while in the appropriate frame these                     
components become Cartan invariants. By choosing this frame basis, it is hoped that the              
algebraic relationships arising from the Bianchi identities, Ricci identities, frame                 
derivatives, or additional requirements will reduce the number of Ricci components                   
involved in the curvature invariants.

\newpage

\section{Basis of scalar curvature invariants in 3D}

\subsection{Algebraically independent scalars}

The number of algebraically independent scalars constructible from the Riemann tensor and its covariant derivatives up to order $p$ is given by \cite{exactSol}:
\beq 
N(n,p) = \frac{n[n+1][(n+{p}+2)!]}{2n!({p}+2)!} - \frac{(n+{p}+ 3)!}{(n-1)!({p}+3)!}+n,\label{thomas}
\eeq
with ${p}\geq 0$, except for $N(2,2)=1$.  

In 3D, we note that $N(3,0) = 3, N(3,1) = 18, N(3,2) = 45, N(3,3) = 87$, $N(3,4) = 147$, and $N(3,5) = 228$.  We remark that
the case of Segre type with no isotropy group (types $\{1,11\}, \{1\bar{z}z\},\{21\},$ and $\{3\}$) thus have $N(3,3) = 87$
algebraically independent curvature invariants.  Similary, those with isotropy groups of dimension $1$, (types
$\{1,(11)\},\{(1,1)1\},$ and $\{(21)\}$ have $N(3,4) = 147$ algebraically independent curvature invariants.  We also note
that for all values of $p$, $N(n,p)$ as given by equation (\ref{thomas}) is equal to the number 
of independent components up to
the $p$th derivative of the Riemann tensor given in \cite{Sousa} minus three (we can deduct three because of the coordinate
conditions that can be imposed).


\subsection{Computing the independent scalar curvature invariants}
As can be seen above, the number of algebraically independent curvature invariants necessary 
is too large to compute by hand. We have therefore utilized the software package 
\emph{Invar}\cite{invar, diffInvar} running on top of Maple. 

Relations among Riemann invariants are a result of symmetries of the Riemann tensor and its covariant derivatives. There exist five such symmetries that are taken into account by \emph{Invar} \cite{diffInvar}:
\begin{enumerate}
\item Permutation Symmetries 
\[R_{bacd}=-R_{abcd},\quad R_{cdab} = R_{abcd}.\]
\item Cyclic Symmetry. 
\[R_{a[bcd]} = 0. \]
\item Bianchi Identity
\[ R_{ab[cd;e]} = 0. \]
\item Commutation of Derivatives. Given any tensor ${T^{a_1\ldots a_n}}_{b_1\ldots b_m}$ we have the identity
\begin{eqnarray*}
&\phantom{}&\nabla_d\nabla_c{T^{a_1\ldots a_n}}_{b_1\ldots b_m} - \nabla_c\nabla_d{T^{a_1\ldots a_n}}_{b_1\ldots b_m} = \\
&\phantom{}&\sum_{k=1}^n {R_{cde}}^{a_k}{T^{a_1\ldots e\ldots a_n}}_{b_1\ldots b_m} - \sum_{k=1}^m {R_{cdb_k}}^e{T^{a_1\ldots a_n}}_{b_1\ldots e\ldots b_m}.\end{eqnarray*}
\item Dimensionally Dependent Identities (also known as Lovelock type identities). Antisymmetrization in $n+1$ indices gives zero (where $n$ is the dimension of the manifold).
\end{enumerate}

Using the commands provided by {\emph{Invar}}, we were able to compute a basis of scalar curvature invariants in
3D.  Let us explicitly present an overdetermined basis of invariants up to $p=1$ using the
cases $[0],[0,0],[0,0,0],$ $[1,1],[1,1,1,1]${\footnote{We refer to the case of invariants
formed by contracting, for example, one undifferentiated Riemann tensor, two once-differentiated
Riemann tensors, and one twice-differentiated Riemann tensor, as $[0,1,1,2]$.  
To give another example,
invariants obtained by contracting two undifferentiated Riemann tensors would be denoted 
by $[0,0]$.}}:

\begin{eqnarray*}
&R, R^{ab}R_{ab}, R^{ab}{R_a}^cR_{bc}, \\
\end{eqnarray*}
\begin{eqnarray*}
&R^{;a}R_{;a}, R^{bc;d}R_{bc;d},R^{bc;d}R_{bd;c},R^{;c}{R_c}^{e;f}R_{;e}R_{;f}, R^{;c}{R_c}^{e;f}{R_{ef}}^{;h}R_{;h}, \\
& R^{;c}{R_c}^{e;f}{{R_e}^{h}}_{;f}R_{;h},R^{;c}{R_c}^{e;f}{{R_f}^h}_{;e}R_{;h}, R^{;c}R^{;f}{R^{hi}}_{;c}R_{hif},\\
& R^{;c}{R_c}^{e;f}{R_e}^{h;i}R_{fh;i},R^{;c}{R_c}^{e;f}{R_e}^{h;i}R_{fi;h},R^{;c}{R_c}^{e;f}{R_e}^{h;i}R_{hi;f},\\
&R^{;c}{R_c}^{e;f}{R_f}^{h;i}R_{hi;e},R^{;c}{R_c}^{e;f}{R^{hi}}_{;e}R_{hi;f},R^{;c}{R^{ef}}_{;c}{R_e}^{h;i}R_{fh;i},\\
&R^{;c}{R^{ef}}_{;c}{R_e}^{h;i}R_{fi;h},R^{;c}{R^{ef}}_{;c}{R_e}^{h;i}R_{hi;f},R^{;c}{R^{ef}}_{;c}{R^{hi}}_{;e}R_{hi;f}, R^{bc;d}{R_{bc}}^{;f}{R_d}^{h;i}R_{fh;i},\\
&R^{bc;d}{R_{bc}}^f{R_d}^{h;i}R_{fi;h},R^{bc;d}{R_{bc}}^{;f}{R_d}^{h;i}R_{hi;f},R^{bc;d}{R_{bd}}^{;f}{R_c}^{h;i}R_{fh;i},R^{bc;d}{R_{bd}}^{;f}{R_c}^{h;i}R_{fi.h},\\
&R^{bc;d}{R_{bd}}^{;f}{R_c}^{h;i}R_{hi;f}, R^{bc;d}{R_{bd}}^{;f}{R_f}^{h;i}R_{hi;c},R^{bc;d}{R_{bc}}^{;f}{R^{hi}}_{;d}R_{hi;f},R^{bc;d}{{R_b}^f}_{;d}{R_c}^{h;i}R_{fh;i},\\
&R^{bc;d}{R_b}^{f;g}{R_{cf}}^{;i}R_{dg;i},R^{bc;d}{R_b}^{f;g}{R_{cg}}^{;i}R_{df;i}, R^{bc;d} {R_b}^{f;g}{R_{cg}}^{;i}R_{di;f}. \\
\end{eqnarray*}

We have included 29 terms involving the first
covariant derivative.  We recall that there are 18 independent components of $R_{ab;c}$.  However, as discussed
earlier, a basis such that all scalar invariants are given as polynomials in the basis contains more
terms.  In table \ref{invarC}, we see the number of independent scalar curvature invariants obtained
per ``case''.  A case in this context refers to the number of Riemann tensors (and their orders of
differentiation) contracted to form scalar curvature invariants.

We note that the total number of curvature invariants returned by \emph{Invar} up to the fourth derivative of the Riemann tensor is actually greater than the number provided by the formula (\ref{thomas}).
\begin{table}
\caption{Number of algebraically independent scalar invariants ($N$) per ``case'' and 
total per derivative ($\bar{N}$). Also see the text and the website \cite{web}. }
\centering
\begin{tabular}{|c|c|c|}
\hline
\phantom{} & \phantom{} & \phantom{} \\
 Case & $N$ & $\bar{N}$ \\
\hline
$[0]$ & 1 &  \\
$[0,0]$ & 1&  \\
$[0,0,0]$& 1& \\
$[0,0,0,0]$& 0& \\
$[0,0,0,0,0]$&0&3 \\
\hline
$[0,0,1,1]$& 15& \\
$[0,0,0,1,1] $&21 & \\
$[1,1]$& 3& \\
$[0,1,1] $&18 & \\
$[1,1,1,1]$& 26&83 \\
\hline
$[0,2]$&2& \\
$[0,0,2] $& 4 & \\
$[0,0,0,2] $&4 & \\
$[0,0,0,0,2] $&4 & \\
$[0,0,2,2]$& 37& \\
$[1,1,2] $&28 & \\
$[0,1,1,2]$& 118 & \\
$[2]$&1& \\
$[2,2]$& 6& \\
$[0,2,2]$&14&  \\
$[2,2,2]$& 30&249 \\
\hline
$[0,1,3]$&18& \\
$[0,0,1,3] $& 50 & \\
$[1,3]$& 4& \\
$[1,2,3] $&96 & \\
$[3,3]$& 8& \\
$[0,3,3]$ & 25&199 \\

\hline
$[0,4]$& 2&  \\
$[0,0,4] $& 5 & \\
$[0,0,0,4] $&7 & \\
$[1,1,4] $&34 & \\
$[2,4]$ & 8 & \\
$[0,2,4]$ & 34 & \\
$[4]$ & 1& \\
$[4,4]$ & 15 &106 \\
\hline
\end{tabular}\label{invarC}
\end{table}

\section{Alternative Bases}  
\subsection{Components}

From equation (\ref{thomas}) and the values of the formula given below equation (\ref{thomas}), we see that there are 3 independent components constructible from $R_{ab}$, $18-3=15$ independent components constructible from $R_{ab;c}$, etc. This means that we can construct a basis consisting of $3$ zeroth order scalar curvature invariants, 15 invariants constructed from $R_{ab;c}$ only, and so on. So a possible basis is
\[ p=0: R, {R_a}^b{R_b}^a, {R_a}^b{R_b}^c{R_c}^a,\]
\[ p=1: R_{ab;c}R^{ab;c}, R_{ab;c}R_{de;f}R^{ab;f}R^{de;c},\ldots,\]
and other independent contractions over indices. We note that such a basis, with higher order (multiple contractions), are extremely difficult to compute.
We could also use table \ref{invarC} to construct such a basis. For example, the cases $[1,1]$ and $[1,1,1,1]$ provide $29$ scalar invariants constructed from $R_{ab;c}$ (alone), the cases $[2],[2,2],[2,2,2]$ provide $37$ scalar invariants constructed from $R_{ab;cd}$, etc.

\subsection{FKWC Basis in 3D}
Field theoretic calculations on curved spacetimes are non-trivial due to the 
systematic occurrence, in the expressions involved, of Riemann 
polynomials. These polynomials are formed from the Riemann tensor by 
covariant differentiation, multiplication and contraction. The results of these 
calculations are complicated because of the non-uniqueness of their 
final forms, since the symmetries of the Riemann tensor as well as the
Bianchi identities can not be used in a uniform manner and monomials 
formed from the Riemann tensor may be linearly dependent in 
non-trivial ways. In~\cite{FKWC1992}, Fulling, King, Wybourne 
and Cummings (FKWC)  
systematically expanded the Riemann polynomials encountered in 
calculations on standard bases constructed from group theoretical 
considerations. They displayed such bases for scalar 
Riemann polynomials of order eight or less in the derivatives of the 
metric tensor and for tensorial Riemann polynomials of order six or 
less. We adopt the FKWC-notations ${\cal R}^r_{s,q}$ 
and ${\cal R}^r_{\lbrace{\lambda_1 \dots \rbrace}}$ to denote, 
respectively, the space of Riemann polynomials of rank r (number of 
free indices), order s (number of differentiations of the metric 
tensor) and degree q (number of factors $\nabla^p 
R_{\dots}^{\dots}$) and the space of Riemann polynomials of rank r 
spanned by contractions of products of the type 
$\nabla^{\lambda_1}R_{\dots}^{\dots}$ 
\cite{FKWC1992}. 
The geometrical identities utilized to 
eliminate ``spurious" Riemann monomials include: 
(i) the commutation of covariant derivatives, 
(ii) the ``symmetry" properties of the Ricci and the 
Riemann tensors (pair symmetry, antisymmetry, cyclic symmetry), and 
(iii) the Bianchi identity and its consequences 
obtained by contraction of index pairs.

\subsubsection{Riemann polynomials of rank 0 (scalars)}

In 3D, the most general expression for a scalar of order six or less in 
derivatives of the metric tensor is obtained by expanding it in the 
FKWC-basis for Riemann polynomials of rank 0 and order 6 or less as follows \cite{FKWC1992}.

{\it The sub-basis for Riemann polynomials of rank 0 and order 2} 
 consists of a single element: $R$ [${\cal R}^0_{2,1}$].

{\it The sub-basis for Riemann polynomials of rank 0 and order 4} 
 has 3 elements: $\Box R$ [${\cal R}^0_{4,1}$]: $R^2$,    
$R_{pq} R^{pq}$,   [${\cal R}^0_{4,2} $].

 {\it The sub-basis for Riemann polynomials of order 6 and rank 
0} in 3D consists of the 10 following elements 
\cite{FKWC1992}: $\Box \Box R$  
$[{\cal R}^0_{6,1}]$: $R\Box R$,   $R_{;p q}R^{pq}$, $R_{pq}\Box R^{pq}$,   
$[{\cal R}^0_{\lbrace{2,0\rbrace}}]$: 
$R_{;p}R^{;p}$,  
$R_{pq;r} R^{pq;r}$,  $R_{pq;r} R^{pr;q}$,
$[{\cal R}^0_{\lbrace{1,1 \rbrace}}]$: $R^3$,   $RR_{pq} R^{pq}$,    
$R_{pq}R^{p}_{\phantom{p} r}R^{qr}$,\\
and so on. These scalar polynomials are not algebraically independent (e.g., $R^2$ is not independent of $R$, $RR_{pq}R^{pq}$ is not independent of $R$ and $R_{pq}R^{pq}$, etc).

\newpage

\section{Example: the Karlhede algorithm}

The maximal order of covariant derivative required for the 
invariant classication of a 3D Lorentzian pseudo-Riemannian manifold ($q \leq 5$ \cite{MW2013}) is relevant
in determining the worst case scenarios for implementing the equivalence algorithm.
It may be possible to provide a minimal basis on a case by case basis; some 
cases may have a smaller $q$ and will then have fewer algebraically independent scalar invariants.
For example,  we could consider Segre type $\{(1,1),1\}$ (Ricci type I or P-type I) with no isotropy group.

In the special case of the (general) P-type I 3D spacetime with invariant
count (3,3), we have that $q=1$. At zeroth order there are 
3 independent scalar invariants (after using the Lorentz freedom). There are
eighteen scalar invariants involving the first covariant derivative, 
of which only 15 are algebraically
independent due to the Bianchi identities \cite{MW2013}. 
Therefore, in order to classify these spacetimes we only need
(after choosing which three terms in the Bianchi identities that are algebraically
dependent), 3+15 = 18 Cartan invariants.

Let us apply the Karlhede equivalence algorithm to the P-type I 3D
spacetimes, in order to identify the Cartan invariants that are
algebraically independent.  We then list some of the simplest scalar
curvature invariants and try to relate them to the algebraically
independent Cartan invariants.  It is argued that, in this
general case, we may use the algebraic independence of the Cartan
invariants to determine the minimal number of algebraically
independent scalar curvature invariants.

\newpage
\subsection{3D spaces of P-type I}	

Following the notation of \cite{MW2013};
at zeroth order, we fix all frame freedom by ensuring $\Psi_0 =
\Psi_4$ and $\Psi_1 = \Psi_3 = 0$, so that $dim H_0 = 0$ - this
ensures that the coframe is invariant; that is, for any diffeomorphism
$\Phi$ we have that $\Phi^* e_a = \tilde{e}_a$, where $e_a$ is the frame basis.  At zeroth order,
there are three non-zero linearly independent components of the Ricci
tensor, namely $R$, $\Psi_0$ and $\Psi_2$.  Hence, at zeroth order
there may be at most three functionally independent invariants
arising.

Let us assume that all three invariants are functionally independent.
Then, in order to complete the algorithm we must set $q=1$ and compute
the covariant derivative of the Ricci tensor.  Instead of working with
the complicated expressions of a rank 3 tensor $R_{ab;c}$ we use the
fact that the frame derivatives of the zeroth order invariants are
already invariants (because the coframe is an invariant coframe) and
simplify the components to produce a set of scalars consisting of
frame derivatives and spin-coefficients:

\beq D \Psi_0, \Delta \Psi_0, \delta \Psi_0, D \Psi_2, \Delta \Psi_2, \delta \Psi_2, D R, \Delta R, \delta R, \alpha, \epsilon, \lambda, \kappa, \pi, \gamma, \tau, \sigma, \nu \nonumber \eeq

\noindent Of these eighteen quantities, only 15 are algebraically independent due to the Bianchi identities \cite{MW2013}. 

This completes the algorithm for the P-type I 3D spacetimes with invariant count (3,3); in order to classify these 
spacetimes we need only choose which three terms in the Bianchi identities are algebraically dependent and list the remaining 3+15 = 18 Cartan invariants. 

\subsection{Cartan Invariants}
Let us consider the minimal basis of algebraically independent scalar curvature invariants.
By completing the Cartan equivalence algorithm we produce two noteable subsets of the Cartan invariants 
by choosing three functionally independent invariants $I_i,~ i\in[1,4]$; we may write the remaining Cartan invariants: 

\begin{itemize} \item $H_{~l}, ~ l \in [1,m]$ : Collection of all invariants that are functionally dependent (f.d.) on $I_i$. 
\item $G_{~l'}, ~ l' \in [1,n],~n \leq m$ : Collection of all invariants that are  algebraically dependent (a.d.) on $I_i$
\end{itemize}

\noindent Treating these as sets we may produce the set of all invariants that cannot be written as algebraic 
expressions of the three invariants $I_i$:

\begin{itemize}
\item $F_{0~l_0},~l_0 \in [1,m-n]$, with $\{F_{0~l_0}\} = \{H_{l}\} \setminus \{G_{l'}\}$. 
\end{itemize}

\noindent Potentially this set may contain invariants that are still algebraically 
dependent on each other, despite being algebraically independent on $I_i$.

To eliminate those elements of $\{F_{0~l_0}\}$ that are algebraically dependent, we repeat the following steps after denoting $I_i = I_i^{(0)}$. 

\begin{enumerate} 
\item Set q=0. 
\item Write $\{F_{0~l_0}\}$ as functions of $I_i^{(0)}$; i.e.,  
$\{F_{0~l_0}\} = \{G_{l_0}(I^{(0)}_i)\}$. Choose four functionally independent invariants from this set 
and denote them as $I_i^{(1)}$. 
\item  Eliminate those invariants which may be expressed algebraically in terms of $I_i^{(1)}$, and denote the subset of invariants that are f.d. but not a.d. on $I_i^{(1)}$ as $F_{1~l_1}$.
\item If $F_{1~l_1}$ has less than three elements, stop the algortihm. Otherwise set q=1 and move to the next step for $q>0$.  
\item Write $\{F_{q~l_q}\}$ as functions of $I_i^{(q-1)}$; i.e.,  $\{F_{q~l_q}\} = \{G_{l_q}(I^{(q-1)}_i)\}$. Choose three functionally independent invariants from this set and denote them as $I_i^{(q)}$ 
\item  Eliminate those invariants which may be expressed algebraically in terms of $I_i^{(q)}$, and denote the subset of invariants that are f.d. but not a.d. on $I_i^{(q)}$ as $F_{q~l_q}$
\item If $F_{q~l_q}$ has less than three elements, stop the algorithm. Otherwise set q=q+1 and repeat steps 5-7  until such a $q$ is found 
\item Given $q_0$ such that $F_{q_0~l_{q_0}}$ has less than three elements, we may express these as functions of $I_i^{(q_0)}$ and determine whether they are algebraically dependent on $I_i^{(q_0)}$ to produce the smaller subset $\tilde{F}_{q_0~l_{q_0}}$ containing at most two invariants which are f.d on $I_i^{(q_0)}$ but not a.d.
\end{enumerate}

\noindent By applying this procedure we produce a list of Cartan
invariants that are algebraically independent:  \beq
F_{\omega~l_\omega} = \displaystyle \bigcup_{j=0}^{q_0} I_i^{(j)}
\bigcup \tilde{F}_{q_0~l_{q_0}} \label{AiCartan} \eeq \noindent The
cardinality of this set determines the number of algebraically
independent polynomial scalar curvature invariants.  Perron's Theorem
\cite{Perron} allows us to reduce the number of algebraically
independent scalar curvature invariants to the cardinality of
$F_{\omega~l_\omega}$, as any scalar curvature invariant may be written
in terms of these quantities; that is, as a polynomial in
$dim(F_{\omega~l_\omega})$-variables.  Of course this result only
applies when there are more polynomials $f_i$ than variables
$x_j,~j\in [1,n]$.

If we wish to determine the smallest number of scalar curvature
invariants $\{ S_i \}$ $i\in [1,N]$ such that any other scalar
curvature invariant may be written as a polynomial in terms of this
set $\{S_i\}$, Perron's theorem is no longer helpful.  As an alternative we
conjecture that it is possible to exploit the relationship between
Cartan invariants and the scalar curvature polynomial invariants in
$\mathcal{I}$-non-degenerate spacetimes, either by directly
expressing the Cartan invariants as polynomials of the scalar
curvature invariants \cite{CH} or indirectly by identifying
algebraically independent polynomials by their dependence on Cartan
invariants.

\subsection{Scalar invariants at zeroth and first order}

As an illustration, we consider one of the `best'-case scenarios for 3D $\mathcal{I}$-non-degenerate spacetimes, the
P-type I spacetimes \cite{MW2013} where the three zeroth-order Cartan invariants are functionally independent.  It should
be emphasized that this is not a proof, but rather an argument towards the feasibility of such an approach.

At zeroth order, the three algebraically independent scalar curvature invariants consist of the Ricci scalar $R = R^a_{~a}$ and 

\beq  S_{ab}S^{ab} = S^2 = 3 \Psi_2^2 + \Psi_0^2,~~and~~S_{ab}S^b_{~d}S^{ad} =  - \Psi_2^3 + \Psi_2 \Psi_0^2 \label{SCI0} \eeq

\noindent where $S_{ab}$ is the trace-free Ricci tensor. 

At first order, we may compute the following scalar curvature invariants of 
lowest degree, $R^1, R^2,R^3$ and $R^4$, respectively:  

\beq R_{,a} R^{,a} = -2 D(R) \Delta(R)+2\delta(R)^2 \nonumber \eeq

\beq \frac14 R_{ab;c}R^{ab;c} &=& 12\tau \Psi_0 \kappa \Psi_2+12\nu \Psi_2 \pi \Psi_0+16 \gamma \Psi_0^2 \epsilon-2\tau\Psi_0^2 \pi-18\nu \Psi_2^2\kappa-12\sigma^2 \Psi_0 \Psi_2+4\sigma \Psi_0^2 \lambda  \nonumber \\
&&+36 \lambda \Psi_2^2 \sigma-12 \lambda^2\Psi_2 \Psi_0-2\kappa \Psi_0^2 \nu-18\pi \Psi_2^2 \tau-3 D(\Psi_2) \Delta(\Psi_2)-\Delta(\Psi_0) D(\Psi_0) \nonumber \\ & &-16\alpha^2\Psi_0^2+\delta(\Psi_0)^2+3\delta(\Psi_2)^2 \nonumber \eeq

\beq \frac12 R_{ab;c}R^{ac;b} &=& -2 \pi \Psi_0 \delta(\Psi_0)-8 \pi \Psi_0^2 \alpha-6 \delta(\Psi_2) \pi \Psi_2-D(\Psi_2) \Delta(\Psi_2) -\Delta(\Psi_0) D(\Psi_0)\nonumber \\
& &-D(\Psi_2) D(\Psi_0)+9 \tau^2 \Psi_2^2 +\nu^2 \Psi_0^2+\kappa^2 \Psi_0^2-\Delta(\Psi_2) \Delta(\Psi_0)+4 \delta(\Psi_2)^2\nonumber \\
& &-18 \nu \Psi_2^2 \kappa-12 \sigma^2 \Psi_0 \Psi_2+4 \sigma \Psi_0^2 \lambda+36 \lambda \Psi_2^2 \sigma -12 \lambda^2 \Psi_2 \Psi_0+4 \Delta(\Psi_2) \lambda \Psi_0\nonumber \\
& &+6 \delta(\Psi_2) \tau \Psi_2-2 \delta(\Psi_2) \nu \Psi_0+4 \Delta(\Psi_2) \gamma \Psi_0-12 \Delta(\Psi_2) \sigma \Psi_2-4 \gamma \Psi_0 D(\Psi_0) \nonumber \\
& &+4 \Delta(\Psi_0) \epsilon \Psi_0-8 \alpha \Psi_0^2 \tau+2 \delta(\Psi_0) \tau \Psi_0-6 \delta(\Psi_0) \nu \Psi_2+6 \kappa \Psi_2 \delta(\Psi_0)-4 D(\Psi_2) \epsilon \Psi_0\nonumber \\
& &+2 \delta(\Psi_2) \kappa \Psi_0 -4 D(\Psi_2) \sigma \Psi_0+12 D(\Psi_2) \lambda \Psi_2+6 \tau \Psi_0 \kappa \Psi_2-2 \tau \Psi_0^2 \pi+16 \gamma \Psi_0^2 \epsilon\nonumber \\
& &+24 \alpha \Psi_0 \nu \Psi_2+24\kappa \Psi_2 \alpha \Psi_0 -6 \tau \Psi_2 \nu \Psi_0+9 \pi^2 \Psi_2^2-6 \kappa \Psi_0 \pi \Psi_2+6 \nu \Psi_2 \pi \Psi_0 \nonumber \eeq


\beq R^a_{~b;a}R^{,b} &=& 4 D(R) \epsilon \Psi_0+ D(R) D(\Psi_0)-2 \delta(R) \kappa \Psi_0+6 \delta(R) \pi \Psi_2+D(\Psi_2) \Delta(R)-2 D(R) \sigma \Psi_0 \nonumber \\
& & +6 D(R) \lambda \Psi_2+4 \delta(\Psi_2) \delta(R)-6 \Delta(R) \sigma \Psi_2+2 \Delta(R) \lambda \Psi_0 +\Delta(\Psi_2) D(R)-6 \delta(R) \tau \Psi_2 \nonumber \\
& & +2 \delta(R) \nu \Psi_0-4 \Delta(R) \gamma \Psi_0+\Delta(R) \Delta(\Psi_0) \nonumber \eeq

These invariants lack the following first order Cartan invariants in each: 

\beq  R_{,a} R^{,a} &:& D \Psi_0, \Delta \Psi_0, \delta \Psi_0, D \Psi_2, \Delta \Psi_2, \delta \Psi_2,\alpha, \epsilon, \lambda, \kappa, \pi, \gamma, \tau, \sigma, \nu \nonumber \nonumber \\
R_{ab;c}R^{ab;c} &:& D(R), \Delta(R), \delta(R).\nonumber \\
R_{ab;c}R^{ac;b} &:& D(R), \Delta(R), \delta(R). \nonumber \\
R^a_{~b;a}R^{,b} &:& R, \delta(\Psi_0), \alpha. \nonumber \eeq

By choosing three non-zero quantities in the Bianchi identities to be algebraically dependent we may remove three more terms
from these three invariants; however, such a choice will depend on the spacetime and the non-vanishing of the Cartan
invariants involved in the Bianchi identities.  

It is clear that the first invariant will be algebraically
independent from the remaining three invariants as the first two lack
all mention of the frame derivatives of $R$, and the third lacks the
quadratic expressions in terms of the frame derivatives of $R$.
Similarly $R^4$ will be algebraically independent to $R^2$ and $R^3$.
Thus by inspection we have found four scalar curvature invariants that
cannot be written as polynomials of each other.  Continuing in this
manner, it is conjectured that this process can be continued for the
list provided from the {\emph{Invar}} package of the algebraically independent
scalar curvature invariants up to first order, and is believed to be the first 18 
scalar curvature invariants listed in subsection (3.2).

\section{Discussion}

In this paper we have discussed a basis for scalar polynomial curvature invariants in 3D
Lorentzian spacetimes.  Results from previous work on
curvature invariants were reviewed.  In particular, we discussed the class of
spaces that can be characterized uniquely by their scalar curvature
invariants, the $\mathcal{I}$-non-degenerate spacetimes \cite{inv}.  Careful consideration
was given to the bound $p$ on the order of the covariant derivatives of the
Riemann tensor that must be included in the contractions that form the
scalar polynomial curvature invariants.  The Karlhede bound
is a well-studied problem in 4D, and some results were reviewed.  Using the
results from \cite{Sousa}, we presented the bound $q = p+1\leq 5$ in 3D.  After
establishing this bound, the computer software \emph{Invar} was used to
calculate an overdetermined basis of scalar curvature invariants in 3D.
Alternative bases, and particularly the so-called FKWC basis,  were discussed in 3D spacetimes. 
In addition, we discussed the equivalence method (using the Karlhede algorithm)
for
computing the Cartan invariants, which can be used to determine a basis of 
all scalar curvature invariants in $\mathcal{I}$-non-degenerate spacetimes.

Future work might include considering 3D spaces of different algebraic types on a case by case basis.  It may be possible to
lower the bounds on $p$ in some cases, as has been the case in 4D by considering spaces of different Petrov types \cite{exactSol}. 
In addition, we could use the equivalence method approach to determine a basis of 
all scalar curvature invariants in other particular 3D spacetimes
\cite{McNutt}. Finally, the
determination of a basis for scalar polynomial curvature invariants in 5D 
would also be an interesting and useful problem.

\newpage

{\em Acknowledgements}. We would like to thank Jose Martin Garcia, Renato Portugal, 
Malcolm MacCallum, and Robert Milson for helpful comments.  
This work was supported, in part, by NSERC of Canada.

\end{document}